\documentclass{amsart}
\usepackage{amsfonts}
\usepackage{verbatim}
\usepackage{amssymb}
\usepackage{amscd}
\usepackage{amsmath}
\usepackage{amsthm}
\usepackage{graphicx}
\usepackage{subfigure}

\newtheorem{theorem}{Theorem}[section]

\newtheorem{conjecture}[theorem]{Conjecture}

\newtheorem{definition}[theorem]{Definition}

\begin{document}
\bibliographystyle{plain}
\title{On Brenti's conjecture about the log-concavity of the chromatic polynomial}
\author{Sukhada Fadnavis}
\email{Sukhada Fadnavis <sukhada@math.stanford.edu>}
\address{Dept. of Mathematics, Stanford University\\ Bldg. 380 \\ Stanford, CA, 94305}
\date{\today}
\maketitle

\begin{abstract}
The chromatic polynomial is a well studied object in graph theory. There are many results and conjectures about the log-concavity of the chromatic polynomial and other polynomials related to it. The location of the roots of these polynomials has also been well studied. One famous result due to A. Sokal and C. Borgs provides a bound on the absolute value of the roots of the chromatic polynomial in terms of the highest degree of the graph. We use this result to prove a modification of a log-concavity conjecture due to F. Brenti. The original conjecture of Brenti was that the chromatic polynomial is log-concave on the natural numbers. This was disproved by Paul Seymour by presenting a counter example. We show that the chromatic polynomial $P_G(q)$ of graph $G$ is in fact log-concave for all $q > C\Delta + 1$ for an explicit constant $C < 10$, where $\Delta$ denotes the highest degree of $G$. We also provide an example which shows that the result is not true for constants $C$ smaller than $1$.
\end{abstract}

\section{Introduction}

Graph coloring is a very well-studied subject. For a graph $G = (V,E)$ we say that a function $\alpha: V \rightarrow \{1, \ldots, q\}$ is a $q$-coloring of $G$ if for each edge $(u,v)$ of $G$ we have $\alpha(u) \neq \alpha(v)$. In general given a graph $G$ it is difficult to say whether it has a $q$-coloring or not, and hence also difficult to count the exactly number of $q$-colorings. Let $P_G(q)$ be the number of $q$-colorings of $G$. If we try to evaluate $P_G(q)$ by inclusion exclusion then we see that $P_G$ is in fact a polynomial known as the \emph{chromatic polynomial}:
\begin{equation}
    P_G(q) = \sum_{E' \subset E} (q)^{C(E')} (-1)^{|E'|},
\end{equation}
where the sum is over all subsets $E'$ of $E$ and $C(E') $ denotes the number of connected components in $E'$. \\

Even though evaluating the chromatic polynomial exactly is a difficult problem in general, many of it's properties have been studied extensively. There is a rich literature about the log-concavity of the chromatic polynomial. For example see \cite{Brenti}, \cite{Brenti2}, \cite{Stanley3}. The roots of the chromatic polynomial have also been extensively studied \cite{Brown}, \cite{BRW}, \cite{SokalZeroes}, \cite{Borgs}, \cite{Fernandez}. \\

Many interesting conjectures about chromatic polynomials can be found in literature. We study here one such conjecture due to F.Brenti. We will fix graph $G$ and use the short notation $P(q)$ to denote the chromatic polynomial whenever there is no ambiguity. In \cite{Brenti} F. Brenti  made the following conjecture about the chromatic polynomial; D. Welsh also made the same conjecture in a private communication with Paul Seymour \cite{Seymour}:

\begin{conjecture}\label{Brenti}(Brenti and Welsh) The chromatic polynomial is log-concave for integer values of $q$ above the chromatic number i.e.\
\begin{equation}
    P(q-1)P(q+1) \leq P(q)^2, \text{ for all } q \geq \chi(G),
\end{equation}
where $\chi(G)$ denotes the chromatic number of $G$.
\end{conjecture}

This conjecture was verified for a large class of graphs in \cite{Brenti}. But it fails to be true in general.

Paul Seymour \cite{Seymour} disproved the conjecture by providing the following counterexample.
\begin{theorem}\label{counterexample} (Paul Seymour \cite{Seymour}) Let $H$ be a graph on $6n$ vertices defined as follows. Consider the vertices partitioned into six equal disjoint subsets $A_1, \ldots, A_6$. For $u \in A_i$ and $v \in A_j$ and $i < j$, there is an edge joining $u$ and $v$ if and only if $(i,j)$ belongs to the set $$\{(1, 2), (1, 3), (2, 3), (2, 4), (3, 4), (1, 5), (3, 5), (1, 6), (2, 6) \}.$$ Then,
\begin{equation}
    \begin{split}
        & P_H(5) \geq 27^n, \\
        & P_H(7) \geq 217^n \text{ and }\\
        & P_H(6) \leq 1080 \times 72^n + 210 \times 64^n + 360\times 48^n + 360\times 36^n + 90\times 16^n.
    \end{split}
\end{equation}

\begin{equation}
    P_H(5)P_H(7) > P_H(6)^2,
\end{equation}
when $n$ is large.
\end{theorem}

Note that in the above example $\chi(H) = 3$ for all $n$. So the original conjecture of Brenti and Welsh places a restriction on the number of colors, $q$, that does not necessarily depend on how large the graph is.

We show that if we change the restriction on $q$ then log-concavity of the chromatic polynomial holds. In particular we show the following result:
\begin{theorem}\label{main}
The chromatic polynomial is log-concave for integer values of $q$ above $C\Delta+1$ i.e.\
\begin{equation}
    P(q-1)P(q+1) \leq P(q)^2, \text{ for all } q \geq C\Delta+1,
\end{equation}
where $\Delta$ denotes the highest degree of $G$, and $C$ is a constant. In particular we know that the above is true for $C = \sqrt{2}K^{*}$ where $K^{*} < 7$ is a constant.
\end{theorem}

We also show that the constant $C$ cannot be smaller than 1 by providing an example in section \ref{example}

Note that in Theorem \ref{counterexample} $H$ has $\chi(H) = 3$ but $\Delta = 5n-5$. Thus, the statement here is much weaker than the original conjecture.

\section{Proof of Theorem \ref{main}}

We shall need the following theorem due to A. Sokal \cite{Sokal}, C. Borgs \cite{Borgs} and Fern\'{a}ndez and Procacci \cite{Fernandez}. It provides a bound on the zeros of chromatic polynomials of general graphs. But before that we need the following definitions:

\begin{definition}Function $F$ is defined as:
    \begin{equation}
     F(a) = \frac{a+e^a}{\log (1+ae^{-a})}.
    \end{equation}
    Define constant $K = \min_a F(a) \leq F(2/5) = 7.964... < 8$.
\end{definition}

\begin{definition}Function $F$ is defined as:
    \begin{equation}
     F(y) = \frac{y}{(2-y)\log y}.
    \end{equation}
    Define constant $K^{*} = \min_{1 < y <2} F(y) = 6.907..$.
\end{definition}

\begin{theorem}\label{SokalZeroes}(\cite{Sokal},\cite{Borgs}, \cite{Fernandez}) Let $G$ be a graph on $n$ vertices with maximum degree $\Delta$. Then,
    \begin{equation}
        |P_G(q)|> 0 \text{ for all $q$ such that } |q| > C\Delta,
    \end{equation}
for a constant $C$.
\end{theorem}
A. Sokal \cite{Sokal} and C.Borgs \cite{Borgs} showed that $C \leq K = 7.964.. <8 $ as above and it was strengthened by Fern\'{a}ndez and Procacci \cite{Fernandez} to show that $C \leq K^{*} = 6.907..<7.$

Now we complete the proof of Theorem \ref{main}

\begin{proof}  Since $P_G(q)$ has real coefficients it can be factored into linear and quadratic real factors. Let's say the real roots of $P$ are $\alpha_1, \ldots \alpha_r$ and the complex roots are $\beta_1, \overline{\beta_1}, \ldots, \beta_s, \overline{\beta_s} $. Hence,
\begin{equation}
    P_G(q) = (q- \alpha_1) \ldots (q-\alpha_r)(q^2 - (\beta_1 + \overline{\beta_1})q + |\beta_1|^2) \ldots (q^2 - (\beta_s + \overline{\beta_s})q + |\beta_s|^2).
\end{equation}

Theorem \ref{SokalZeroes} gives us the bounds,
\begin{equation}
    |\alpha_i| \leq K^{*}\Delta \text{ and } |\beta_j| \leq K^{*}\Delta,
\end{equation}
for $1 \leq i \leq r$ and $1 \leq j \leq s$.

Note that,
\begin{equation}
    (q-1- \alpha_i)(q +1 -\alpha_i) =  (q-\alpha_i)^2 -1 <  (q-\alpha_i)^2,
\end{equation}
 and also both the LHS and RHS are positive since $|\alpha_i| \leq K^{*}\Delta \leq q - 1$ when $q > \sqrt{2}K^{*}\Delta +1$.

 Now let $p(q) = q^2 - (\beta_j +\overline{\beta_j})q + |\beta_j|^2.$ Note that $p(q) =  |q - \beta_j|^2 > 0$ for all real $q$. Hence it suffices to prove that,
 \begin{equation}
    p(q)^2 \leq p(q-1)p(q+1) \text{ for } q > \sqrt{2}K^{*}\Delta +1.
 \end{equation}

To see this let $\beta_j = a_j + ib_j$ for $a_j, b_j$ real. Then,
$$q^2 - (\beta_j + \overline{\beta_j})q + |\beta_j|^2 = q^2 - (2a_j)q + (a_j^2 + b_j^2) = (q-a_j)^2 + b_j^2.$$

Now,
\begin{equation}
    \begin{split}
        & p(q-1)p(q+1) \\
        & = ((q-1-a_j)^2 + b_j^2)((q+1-a_j)^2 + b_j^2)\\
        & = ((q-a_j)^2 + b_j^2 + 1 - 2q + 2a_j)((q -a_j)^2 + b_j^2 + 1 + 2q - 2a_j) \\
        & = ((q-a_j)^2 + b_j^2)^2 + 2((q-a_j)^2 + b_j^2) + 1 - 4(q -a_j)^2  \\
        & = ((q-a_j)^2 + b_j^2)^2 + 1 -2(q-a_j)^2 + 2b_j^2 \\
        & \leq ((q-a_j)^2 + b_j^2)^2 + 1 -2(q-a_j)^2 + 2((K^{*}\Delta)^2 - a_j^2) \\
        & = ((q-a_j)^2 + b_j^2)^2 + 1  + 2(K^*\Delta)^2 - 2(a_j^2 + (q-a_j)^2) \\
        & \leq ((q-a_j)^2 + b_j^2)^2 + 1 + 2(K^*\Delta)^2 - 4(q/2)^2 \\
        & \leq ((q-a_j)^2 + b_j^2)^2 \text{ since } |q| > \sqrt{2}K^{*}\Delta + 1 \\
        & = p(q)^2.
    \end{split}
\end{equation}

The first inequality above is true because $ |\beta_j|^2 = a_j^2 + b_j^2 \leq 2(K^*\Delta)^2$. The second inequality is true since $2a^2 + 2b^2 \geq (a+b)^2 $ for real numbers $a,b$.

Multiplying together the above linear and quadratic inequalities we get the desired result.
\end{proof}

\section{Example}\label{example}

In this section we show that the constant in theorem \ref{main} cannot be less that 1. To show this we construct a slight modification of the counterexample due to Paul Seymour \ref{counterexample}. 

Consider a graph $S$ on $n^2$ vertices as follows. Let $A_1, \ldots, A_6$ be disjoint sets of vertices such that $|A_i| = n$ for all $i$. For $u \in A_i$ and $v \in A_j$ and $i < j$, there is an edge joining $u$ and $v$ if and only if $(i,j)$ belongs to the set $$\{(1, 2), (1, 3), (2, 3), (2, 4), (3, 4), (1, 5), (3, 5), (1, 6), (2, 6) \}.$$ Further, let $v_1, \ldots, v_{n^2-6n}$ be the remaining $n^2-6n$ vertices with edges from $v_i$ to all other vertices in the graph, for all $i$. Thus, maximum degree of the graph is $\Delta = n^2-1$. Note that,
\begin{equation}
    P_S(k) = \binom{k}{n^2-6n} \times (n^2-6n)! \times P_H(k-n^2-6n),
\end{equation}
where $H$ is the graph in \ref{counterexample}. This can be seen as follows: choose colors for vertices $v_1, \ldots, v_{n^2-6n}$. The remaining vertices induce a graph isomorphic to $H$ and it needs to be colored with the remaining $k-n^2-6n$ colors. Hence by \ref{counterexample} we have,
\begin{equation}
    \begin{split}
        &P_H(n^2-6n+5) \geq \binom{n^2-6n+5}{n^2-6n} \times (n^2-6n)! \times 27^n, \\
        &P_H(n^2-6n+7) \geq \binom{n^2-6n+7}{n^2-6n} \times (n^2-6n)! \times 217^n \text{ and }\\
        &P_H(n^2-6n+6) \leq \\
        &\binom{n^2-6n+6}{n^2-6n} \times (n^2-6n)!\times (1080 \times 72^n + 210 \times 64^n + 360\times 48^n + 360\times 36^n + 90\times 16^n).
    \end{split}
\end{equation}
Thus for large $n$ one has,
\begin{equation}
    P_H(n^2-6n+6)^2 \leq P_H(n^2-6n+5)P_H(n^2-6n+7). 
\end{equation}
Since in this example $\Delta = n^2-1$ and,
\begin{equation}
    \lim_{n\rightarrow \infty} \frac{n^2-6n+7}{n^2-1} = 1,
\end{equation} 
we see that the constant in theorem \ref{main} cannot be smaller than 1.

\section{Further comments}

The example above shows that if the constraint on $q$ is stated in terms of $\Delta$ then $C$ cannot be smaller than 1. But it is still not clear if the correct constraint on $q$ should be stated in terms of $\Delta$. In the above example $\chi(S) = n^2-3n$. So the possibility that $P_G(q)$ is log-concave for $q>C\chi(G)$ is not ruled out. This result would be stronger than the result we prove here since $\chi(G) < \Delta +1$ and many times much smaller. We think that such a result will not be true. It will be interesting to find a family of graphs with bounded chromatic numbers but for which the chromatic polynomials fail to be log-concave for bigger and bigger values of $q$.

\section{Acknowledgements}
The author would like to thank Alan Sokal for pointing out the work of Fernandez and Procacci \cite{Fernandez}and other useful comments.

\bibliography{BrentiConjecture}{}

\end{document}